\documentclass[10pt]{article}
\usepackage{amsmath}
\usepackage{latexsym}
\usepackage{amsfonts}
\usepackage{amssymb}

\newcommand{\q}[1]{\mbox{$#1$}}
\newcommand{\op}{\mbox{$\omega$}}

\newcommand{\B}[1]{\mbox{\boldmath $#1$}}

\newcommand{\g}[1]{\mbox{$\; \ulcorner #1  \urcorner$\,}}
\newcommand{\ol}[1]{\mbox{$\overline{#1}$}}
\newcommand{\pr}{\mbox{$\vdash_{\mathcal C}$}} 

\newcommand{\xy}[2]{\mbox{$ #1 \, \times \, #2$}}

\newcommand{\ppa}[2]{\parbox[c]{40pt}{\nhp\begin{picture}(40,40)
\put(5,16){\vector(1,0){35}} \put(5,24){\vector(1,0){35}}
\put(20,7){$\scriptstyle #1$} \put(20,29){$\scriptstyle #2$}
\end{picture}}\hme }

\newcommand{\etri}[6]{
\begin{picture}(120,120)
\thicklines \setlength{\unitlength}{.6\unitlength}
\put(15,162){#1} \put(47,167){\vector(1,0){95}} \put(162,162){#2}
\put(30,90){#6} \put(150,90){#5} \put(85,20){#3} \put(90,175){#4}
\put(27,150){\vector(1,-2){55}} \put(160,150){\vector(-1,-2){55}}
\end{picture} }

\newcommand{\nf}{\hspace*{-5pt}}

\def \y {\'{\i}}
\def \cao {\c c\~ao}

\def \coes {\c c\~oes}

\def \ao {\~ao }

\def \beq { \begin{equation} }
\def \eeq { \end{equation} }

\def\overset#1\to#2{\mathrel{\mathop{#2}\limits^{#1}}}
\def\underset#1\to#2{\mathrel{\mathop{#2}\limits_{#1}}}

\def \varinjlim {\underset{\longrightarrow}\to{lim}}
\def \varprojlim {\underset{\longleftarrow}\to{lim}}

\def \li {\varinjlim}
\def \lp {\varprojlim}

\def \lli#1 {\mathrel{\mathop{\li}\limits_{#1}}}
\def \llp#1 {\mathrel{\mathop{\lp}\limits_{#1}}}

\def \empty {\emptyset}

\newcommand{\hem}{\hspace*{1em}}
\newcommand{\hme}{\hspace*{.5em}}
\newcommand{\hmm}{\hspace*{2em}}
\newcommand{\vso}{\medskip}
\newcommand{\vsm}{\medskip}
\newcommand{\vst}{\bigskip}

\newcommand{\hfl}{\hspace*{\fill}}

\newlength{\dede}     

\newcommand{\hp}{\hspace*{\parindent}} 
\newcommand{\nhp}{\hspace*{-\parindent}}
\newcommand{\po}[2]{\mbox{$#2^{#1}$}}
\newcommand{\tf}{\mbox{$f$}}

\newcommand{\dt}{\mbox{$\cdot$}} 

\newcommand{\lra}{\mbox{$\longrightarrow$}}  
\newcommand{\com}{\mbox{$\circ$}} 

\newcommand{\rs}[2]{\mbox{$ #1_{| #2}$}}  
\newcommand{\se}[2]{\mbox{$ #1_{#2}$}}

\newcommand{\equ}{\mbox{$\equiv$}}

\newcommand{\en}{\mbox{$\in$}}
\newcommand{\nen}{\mbox{$\not \in$}}

\newcommand{\Int}{\mbox{$\displaystyle \cap$}}
\newcommand{\niq}{\mbox{$\not =$}}
\newcommand{\npl}[1]{\mbox{$ #1_{1}, \ldots , #1_{n}$}}  

\newcommand{\prd}{\mbox{$\prod$}}  
\newcommand{\cpd}{\mbox{$ \coprod$}} 

\newcommand{\sub}{\mbox{$\subseteq$}}

\newcommand{\bcup}{\mbox{$ \bigcup$}}

\newcommand{\bcap}{\mbox{$ \bigcap$}}

\def \empty {\emptyset}
\newcommand{\0}{\mbox{$\emptyset$}} 

\newcommand{\ga}{\mbox{$\gamma$}}

\newcommand{\Si}{\mbox{$\Sigma$}}

\newcommand{\fh}{\raisebox{.4ex}{\mbox{$\varphi$}}}

\newcommand{\lam}{\mbox{$\lambda$}}

\newcommand{\de}{\mbox{$\delta$}}

\newcommand{\al}{\mbox{$\alpha$}}
\newcommand{\bt}{\mbox{$\beta$}}

\newcommand{\tet}{\mbox{$\theta$}}

\newcommand{\ta}{\mbox{$\tau$}}
\newcommand{\et}{\mbox{$\eta$}}
\newcommand{\ps}{\raisebox{.4ex}{\mbox{$\psi$}}}

\newcommand{\cC}{\mbox{$\cal C$}}

\newcommand{\cF}{\mbox{$\cal F$}}
\newcommand{\cI}{\mbox{$\cal I$}}

\newcommand{\cA}{\mbox{$\cal A$}}

\newcommand{\cM}{\mbox{$\cal M$}}

\newcommand{\cU}{\mbox{$\cal U$}}

\newcommand{\ex}{\mbox{$\exists$}}
\newcommand{\fa}{\mbox{$\forall$}}
\newcommand{\all}{\mbox{$\forall$}}
\newcommand{\w}{\mbox{$\wedge$}} 
\newcommand{\jo}{\mbox{$\vee$}} 

\newcommand{\rra}{\mbox{$\rightarrow$}}  

\newcommand{\mo}{\mbox{$\models$}}

\newcommand{\rla}{\mbox{$\leftrightarrow$}}
\newcommand{\Ra}{\mbox{$\Rightarrow$}}
\newcommand{\Lra}{\mbox{$\Leftrightarrow$}}

\def \ra {\rightarrow}

\def\overset#1\to#2{\mathrel{\mathop{#2}\limits^{#1}}}
\def\underset#1\to#2{\mathrel{\mathop{#2}\limits_{#1}}}

\newcommand{\srel}{\stackrel} 
\newcommand{\und}[1]{\raisebox{-.2ex}{\underline{\raisebox{.2ex}{#1}}}}  

\newcommand{\m}{\mbox{$\leq$}} 
\newcommand{\G}{\mbox{$\geq$}} 

\newcommand{\up}[1]{\mbox{$ #1^{\rightarrow}$}}
\newcommand{\dn}[1]{\mbox{$ #1^{\leftarrow}$}}

\newcommand{\fm}[1]{\mbox{$\langle \, #1 \,  \rangle$}}

\newcommand{\ili}{\raisebox{-2pt}{$\srel{\textstyle{lim}}{{}_{{}_{{}_{\longrightarrow}}}}$}} 
\newcommand{\dli}{\raisebox{-2pt}{$\srel{\textstyle{lim}}{{}_{{}_{{}_{\longleftarrow}}}}$}} 

\newcommand{\fo}{\mbox{$\bf 1\!\!1$}}    

\def \varinjlim {\underset{\longrightarrow}\to{lim}}
\def \varprojlim {\underset{\longleftarrow}\to{lim}}

\def \li {\varinjlim}
\def \lp {\varprojlim}

\def \lli#1 {\mathrel{\mathop{\li}\limits_{#1}}}
\def \llp#1 {\mathrel{\mathop{\lp}\limits_{#1}}}

\newtheorem{Th}{Theorem}
\newtheorem{Co}[Th]{Corollary}
\newtheorem{Df}[Th]{Definition}
\newtheorem{Pro}[Th]{Proposition}
\newtheorem{Le}[Th]{Lemma}
\newtheorem{Exa}[Th]{Example}
\newtheorem{Rem}[Th]{Remark}
\newtheorem{Fa}[Th]{Fact}
\newtheorem{Que}[Th]{Question}
\newtheorem{Ct}[Th]{}
\newtheorem{Afi}[Th]{Claim}

\newcommand{\baf}{\begin{Afi}\nf{\sl }}
\newcommand{\eaf}{\end{Afi}}

\newcommand{\bdf}{\begin{Df}\nf{\bf }}
\newcommand{\edf}{\end{Df}}
\newcommand{\bte}{\begin{Th}\nf{\bf }}
\newcommand{\ete}{\end{Th}}  
\newcommand{\bco}{\begin{Co}\nf{\bf }}
\newcommand{\eco}{\end{Co}}
\newcommand{\ble}{\begin{Le}\nf{\bf }}
\newcommand{\ele}{\end{Le}}
\newcommand{\bpr}{\begin{Pro}\nf{\bf }}
\newcommand{\epr}{\end{Pro}}
\newcommand{\bex}{\begin{Exa}\nf{\bf } \rm}
\newcommand{\eex}{\end{Exa}}
\newcommand{\bre}{\begin{Rem}\nf{\bf } \rm}
\newcommand{\ere}{\end{Rem}}
\newcommand{\bfa}{\begin{Fa}\nf{\bf } \sl}
\newcommand{\efa}{\end{Fa}}
\newcommand{\bqt}{\begin{Que}\nf{\bf }}
\newcommand{\eqt}{\end{Que}}
\newcommand{\bdm}{\nhp{\bf Proof.  }}

\newcommand{\qdr}{\hfl $\square$ }
\newcommand{\bct}{\begin{Ct}\nf{\bf } \rm}
\newcommand{\ect}{\end{Ct}}
\newcommand{\bxa}{\begin{Exa}\nf{\bf } \rm}
\newcommand{\exa}{\end{Exa}}

\makeindex

\begin{document}

\title{Profinite Structures are Retracts of\\
 Ultraproducts of Finite Structures}

\author{ {\em H. L. Mariano} --- Centro de L\'ogica, Epistemologia
e Hist\'oria da Ci\^encia, \\
Universidade de Campinas, S\ao\ Paulo, S.P., Brazil. \\
E-mails: hugomar@cle.unicamp.br, hugomar@ime.usp.br}

\date{December 2003}

\maketitle

\section*{Abstract}

\bigskip

  We establish the following model-theoretic characterization: profinite $L$-structures, the
cofiltered limits of finite $L$-structures, are retracts of
ultraproducts of finite $L$-structures. As a consequence, any
elementary class of $L$-structures axiomatized by $L$-sentences of
the form $\forall \vec{x} (\psi_{0}(\vec{x}) \ra
\psi_{1}(\vec{x}))$, where $\psi_{0}(\vec{x}),\psi_{1}(\vec{x})$
are existencial-positives $L$-formulas, is closed under the
formation of profinite objects in the category {\bf L-mod}, the
category of structures suitable for the language $L$ and
$L$-homomorphisms.

\section{Introduction}

\hp The results presented here belong to the interface between
Category Theory and Model Theory. These results are contained in
{\em Chapter 2} of \cite{Mrn1}. Our primary motivation was
\cite{KMS}, a paper which introduces the class of direct limits of
finite abstract order spaces. The Theory of Spaces of Orderings is
an axiomatization of the algebraic theory of quadratic forms on
fields (see \cite{Mar1}). Later, in \cite{DM2}, it was presented a
{\em first-order} axiomatization of the algebraic theory of
quadratic forms, the Special Groups Theory, which is, in some
sense, a dual approach to the Theory of Orderings Spaces, but with
an advantage: it permits an approach of quadratic forms theory by
the logical methods of Model Theory.

Detailing the work:

We consider $L$, a first-order language with equality. We denote
{\bf L-mod}, the category of all structures suitable to the
language $L$ and $L$-homomorphisms. As preparation we present some
species of limits and  colimits in the category {\bf L-mod} and we
relate one of the principals constructions in Model Theory, the
notion of reduced product of structures, with the categorial
constructions of product and filtered colimit in {\bf L-mod}
(Proposition \ref{colu-pr}). Our main result, the Theorem
\ref{repro-te}, claims that the {\em profinite} $L$-structures,
the cofiltereds limits of finite $L$-structures, are {\em retracts
of ultraproducts} of finite $L$-structures. As a consequence, each
elementary class of $L$-structures axiomatized by $L$-sentences
like $\all \vec{x} (\psi_{0}(\vec{x}) \ra \psi_{1}(\vec{x}))$,
were $\psi_{0}(\vec{x}),\psi_{1}(\vec{x})$ are
existencial-positives $L$-formulas, is closed under the formation
of profinite objects in the category {\bf L-mod} (Corollary
\ref{repro-co}).

Applying the central results, \ref{repro-te} and \ref{repro-co},
to the Special Groups Theory we conclude that there are profinite
special groups and that they are retracts of ultraproducts of
finite special groups (Section \ref{fin-sec}).

\section{Preliminaries}

\hp We assume familiarity with the basic notions of Category
Theory (category, functor, natural transformation,
limits/colimits, ...) and of Model Theory (language, structure,
homomorphism, elementary embedding, reduced products, ...). Our
reference about Category Theory is \cite{Mac}; for Model Theory we
use \cite{CK} e \cite{BS}.

We clarify below some topics needed to the development of the
results obtained in this work.

\subsection{Retracts} \label{retr-ssec}

\hp Let \cC\ be a category and $A, B$ objects of ${\cC}$. $A$ is
called a {\em retract} of $B$ when there are morphisms $s : A \ra
B$ and $r : B \ra A$ such that $r \circ s = Id_{A} : A \ra A$. In
this case we say that $r$ is a retraction and $s$ a section.

Is immediate to verify that any section is a monomorphism and,
dually, any retraction is an epimorphism; that a morphism is
invertible (or isomorphism) precisely when it is  simultaneously a
section and a retraction.

We remark that the proposition: ``all epimorphism in the category
of sets and functions ({\bf Set}) is a retraction'' is equivalent
to the Axiom of Choice.

\subsection{Directed Sets}  \label{conjdir-ssec}

\hp Let \fm{I, \m}\ be a poset, i.e. \m \ $\sub$ \ $I \times I$ is
a binary relation that is reflexive, symmetric and transitive in
the set $I$. For each \q i \en\ \q I we define \ \dn i  =  \{\q j
\en\ \q I : \q j \m\ \q i\} \ , \ \up i  =  \{\q j \en\ \q I : \q
i \m\ \q j\}. \
We say that: \\
$\ast$ \ \fm{I, \m}\  is {\em upward directed } (or {\em
filtered}) if $I \neq\empty$ and for each \q i,
\q j \en\ \q I, \, \up i \Int\ \up j \niq\ \0. \\
$\ast$ \ \fm{I, \m}\ is {\em downward directed } (or {\em
cofiltered}) if $I \neq\empty$ and for each \q i, \q j \en\ \q I,
\, \dn i \Int\ \dn j \niq\ \0.

Clearly a poset \fm{I, \m}\ is upward directed iff its opposite
poset $\fm{I,\m }^{op}$\ is downward directed and vice-versa. When
we make mention to directed posets we {\em always} will be
refiring the upward directed orders.

We say that a filter ${\cF}$ in the set \q I is a {\em directed
filter} in the poset \fm{I, \m}\ when, for each $i \in I$, we have
${i}^{\ra} \in {\cF}$.

\ble \label{ufildir-le} If \fm{I,\m}\ is a directed poset then
there is a  directed \und{ultrafilter} ${\cU}$ in \fm{I,\m}.

\ele

\bdm Because \fm{I, \m}\ is directed we verify, by induction on
$\q n \in \omega$ , that for each $\{ i_{0},\ldots, i_{n-1} \} \
\sub\ $ \q I \ there is $j \in I$ such that \, $j \, \G\ i_{0} ,
\ldots i_{n-1}$ , so \ $\empty \ \neq \ \up{j} \  \sub \
\bcap_{m<n}$ \up{i_m} \ . Hence the set  \q S \ = \ \{\up i: \q i
\en\ \q I\} \ has the finite intersection property and then there
is an ultrafilter {\cU} \ such that \q S \sub\ {\cU}. \qdr

\subsection{Pure Morphisms} \label{Lpur-ssec}

\hp Let \q L an arbitrary first-order language with equality.

\bdf \label{sent-df} A formula \fh\ in the language \q L is called: \\
$\ast$ {\em positive} if the symbols of implication and negation
do not occur in \fh; \\
$\ast$ {\em existencial positive (e.p.)} if it is obtained from
the atomic formulas by the connectives \w, \jo\ and the
existencial quantifier \ex; \\
$\ast$ {\em positive primitive (p.p.)} if it is written like
\ex\,\ol x\,\fh, where \fh\ is a conjunction of atomic formulas.

\edf

We denote: \\
$\ast$ \B{\po{+}{\ex}(L)} the set of all \q L-formulas that are
logically equivalents, in the classical predicate calculus, to a
formula
existencial positive in \q L; \\
$\ast$ {\bf pp(\B L)} the set of all \q L-formulas that are
logically equivalents, in the classical predicate calculus, to a
formula positive primitive in \q L.

By induction on complexity, we see that if \fh\ \en\
$\po{+}{\ex}(L)$ then there are finite subsets \npl P of
${\rm pp}(L)$ such that \\
\hfl  \pr\ \fh\ \, \rla\ \, $(\ps_1\ \jo\ \ps_2\ \jo\ \ldots \jo\
\ps_n)$,  \hfl \\
where \se{\ps}{j}\ is a conjunction of formulas in \se P j, \q 1
\m\ \q j \m\ \q n.

\bdf \label{Lpur-df} A function between $L$-structures, \mbox{\q f
: \q M \lra\ \q N}, is called a {\em pure $L$-morphism} if for
each formula $\fh(\npl v)$ \en\ $\po{+}{\ex}(L)$ and \ol a \en\
\po n M

\nhp \hfl \q M \mo\ $\fh[\ol a]$ \, \Lra\ \, \q N \mo\ $\fh[f\ol
a]$. \hfl

\edf

It is not difficult to see that: a function between $L$-structures
is a pure $L$-morphism iff it is a $L$-homomorphism that reflects
the validity of formulas in $pp(L)$; all pure $L$-morphism is a
$L$-imbedding; all elementary $L$-imbedding and all $L$-section
\footnote{A $L$-section is a $L$-homomorphism that admits a
retraction that is also a $L$-homomorphism.} is a pure
$L$-imbedding. We register also the following

\ble \label{Lpur-le} Let \Si\ be a set of \q L-sentences of the
form $\all \vec{x} (\psi_{0}(\vec{x}) \ra \psi_{1}(\vec{x}))$,
where $\psi_{0}(\vec{x})$, $\psi_{1}(\vec{x})$ \en\
$\po{+}{\ex}(L)$. Let \q N be a \q L-structure such that $N
\models \Sigma$; if \q M is a \q L-structure and there is a pure
\q L-morphism from \q M to \q N then $M \models \Sigma$. \qdr \ele

\section{The category L-mod}

\hp Henceforth we fix \q L an arbitrary first-order language with
equality. We shall write $ct(L)$ for the set of all symbols for
constants of the language and, for each \q n \G\ 1, $op(n,L)$
denotes the set of all symbols for operations with aridity n and
$rel(n,L)$ for the set of all symbols for n-ary relations.

We denote {\bf L-mod} the {\em category} of all structures
suitable to the language \q L and of \q L-homomorphisms between
them\footnote{We will {\em not exclude} here the possibility of a
\q L-structure be empty. As a structure is non-empty iff it
satisfies the sentence $\exists v_{0} (v_{0} = v_{0})$ we should
write the instantiation axiom as \ \fa\q v \fh\ $\wedge\ \exists
v_{0} (v_{0} = v_{0})$  \rra\ \, \fh(\g{\ta \mid v}) \ were \ta\
is a term free for \q v in \fh  .}.

{\bf L-mod} is a complete and cocomplete category, i. e., all
diagram $D : {\cI} \  \lra \ {\bf L-mod}$, where ${\cI}$ is a
small category, is base of some limit cone and some colimit
co-cone (\cite{Mac}, Chapter 5).

We will detail below some of that categorials constructions and
how the reduced products, one of the fundamentals notions of Model
Theory, is related with these constructions.

\subsection{Limits in L-mod} \label{Llim-ssec}

\bct {\bf Products in L-mod:}  \label{Lprod-ct} Let $I$ a set and
$\{ {M}_{i} \ : i \in I \}$, a family of \q L-structures. We
consider $M\ =\ \prd_{i \in I}\ \se{M}{i}$ the product of their
underlying sets. We make \q M a \q L-structure, defining the \q
L-symbols interpretations coordinate-wise. More explicitly, for
each natural \q n \G\ 1 :

\nhp $\ast$ If \q c \en\ $ct(L)$, \po{M}{c}\ \, = \,
$\fm{\po{M_{i}}{c}}_{i \in I}$;

\nhp $\ast$ Se \op\ \en\ $op(n,L)$ e \fm{\npl{s}}\ \en\ \po{n}{M},
then

\nhp \hfl \po{M}{\op}(\npl{s})  \, = \, $\fm{\po{M_{i}}{\op}(\se s
1(\q i), \se s 2(\q i), \ldots, \se s n(\q i))}_{i \in I}$; \hfl

\nhp $\ast$ If \q R \en\ $rel(n,L)$ e \ol s  \en\ \po{n}{M}, then

\nhp \hfl \q M \mo\ \q R[\ol s] \hem \Lra\ \hem \fa\ \q i \en\ \q
I, \, \se M i \mo\ \q R[\se s 1(\q i), \se{s}{2}(\q i), \ldots,
\se{s}{n}(\q i)]. \hfl

By induction on the complexity of terms and formulas we get:

\nhp (A) If \ta(\npl v) is a term em \q L e \ol s \en\ \po n M,
then

\nhp \hfl \po{M}{\ta}(\ol s) \, = \, $\fm{\po{M_i}{\ta}(\se s
1(i), \ldots, \se s n(i))}_{i \in I}$. \hfl

\nhp (B) If \fh(\npl v) is a atomic formula in \q L and \ol x \en\
\po n M, then

\nhp \hfl \q M \mo\ \fh[\ol s] \hem \Lra\ \hem \fa\ \q i \en\ \q
I, \, \se M i \mo\ \fh[\se s 1(\q i), \ldots,  \se{s}{n}(\q i)].
\hfl

 Observe that the canonicals projections, \se{\pi}{i}\ : \q M \lra\
\se{M}{i}, $i \in I$, are \q L-morphisms. It is easy to see that
this construction is the product of the family $\{ \se M i : i \in
I \}$ in the category {\bf L-mod}.

Particularly, when $I = \empty$, we have the

\nhp \und{Final object of {\bf L-mod} :}
Let \fo\ = \{\0\}, where all \q n-ary relation symbols are
interpreted by \po{n}{\fo}, all \q n-ary functional symbols are
interpreted as the unique function $\fo^{n} \ra \fo$ and all
constant symbols are interpreted as the unique element of \fo\ .
We see, by induction on the complexity, that all $L$-formula
positive (Definition \ref{sent-df}) is satisfiable em \fo\ ; hence
all $L$-sentence of the form $\all \vec{x} (\psi_{0}(\vec{x}) \ra
\psi_{1}(\vec{x}))$, where $\psi_{0}(\vec{x}),\psi_{1}(\vec{x})$
are positive $L$-formulas, is true in \fo\ .

We observe also that $ \prd_{i \in I}\ \se{M}{i} \ = \ \empty$ \
iff \ there is a $i \in I$ such that $M_{i} = \empty$ .\footnote{A
version of choice axiom.}

\ect

\bct {\bf Equalizers in L-mod:} \label{Lequa-ct} Let \q D = (\q A
\ppa{g}{f}\ \q B) be \q L-morphisms. We define

\nhp \hfl \q E = \{\q a \en\ \q A : $f(a)$ = $g(a)$\}. \hfl

\nhp If \q c \en\ $ct(L)$,  \po{E}{c}\ $=_{def}$ \po{A}{c}\ \en\
\q E. Further, if \op\ \en\ $op(n, L)$ and \ol{a}\ \en\ \po{n}{E},
then

\nhp \hfl \q f(\po{A}{\op}(\ol{a})) = $\po{B}{\op}(\q f(\ol{a}))$
= $\po{B}{\op}(\q g(\ol a))$ = \q g(\po{A}{\op}(\ol a)), \hfl

\nhp and \q E is closed with respect to functional symbols
interpretations. For each \q R \en\ $rel(n,L)$, let
\mbox{\po{E}{R}\ = \q R \Int\ \po{n}{E}}. So the canonical
inclusion, \et\ : \q E \lra\ \q A, is a \q L-imbedding. Moreover,
because $f\ \com\ \et\ =\ g\ \com\ \et$, $(E; \{\et, \, f\ \com\
\et\})$ is a cone over \q D em {\bf L-mod}. This cone is the
equalizer of (\q f, \q g).

\ect

From the remarks \ref{Lprod-ct} and \ref{Lequa-ct} below and the
construction of limits from products and equalizers (\cite{Mac},
section 5.2) we get the

\bco \label{Llim-co} Let \  $D  \ : \ {\cI} \ \lra {\bf L-mod}$ \
: \ $(i \overset{\alpha}\to\lra j)$  \ $\mapsto$ \ $(D_i
\overset{f_\alpha}\to\lra D_j)$  \ be a {\cI}-diagram in {\bf
L-mod} and \ $( \ M \overset{\lam_{i}}\to\lra D_i \ : \ i \in
Obj({\cI}) \ ) $ \ be a cone over $D$. We take \lam\ = \ $(
\se{\lam}{i})_{i \in I}$ \ the unique function from \q M to
$\prd_{i \in I} D_i$, such that \se{\pi}{i}\ \com\ \lam\ =
\se{\lam}{i}, for each \q i \en\ $Obj(\cI)$ . Then \q M \ is
$(\text{isomorphic to})$ \dli\ \q D \ iff

\vsm \nhp $[lim\ 1]$ : The image of \lam\ in $\prd_{i \in I} D_i$
\ is the set

\vso \nhp \hfl  \{\q x \en\ $\prd_{i \in I} D_i$ : for all arrow
of ${\cI}$ , $(i \overset{\alpha}\to\lra j)$ , we have \
$\se{\tf}{\alpha}(\se{\pi}{i}(x))$ = $\se{\pi}{j}(x)$\}. \hfl

\vsm \nhp $[lim\ 2]$ : If $\fh(\npl{v})$ is an atomic formula in
\q L and  \ol s \en\ \po{n}{\q M},

\vso \nhp \hfl  \q M \mo\ $\fh[\ol s]$ \, \Lra\ \, \fa\ \q i \en\
$Obj({\cI})$ , \, $D_i$ \mo\ $\fh[\se{\lam}{i}(\ol s)]$. \hfl \qdr
\eco

\subsection{Reduced Products and Ultraproducts of L-structures}
\label{redprod-ssec}

\bct \label{redprod-ct} Let \q I be a {\em non empty} set and
\{\se M i : \q i \en\ \q I\} be a family of \q L-structures, {\em
all non empty}. We fix \cF\ a filter in \q I and we consider \q M
= $\prd_{i \in I}$  \se M i the product of their underlying sets
(so $M \neq \empty$).
 We define a binary relation \se{\tet}{F}\ in
\q M :

\nhp \hfl \q x \se{\tet}{F}\ \q y \hem \Lra\ \hem \{\q i \en\ \q I
: $x(i)$ = $y(i)$\} \, \en\ \, \cF . \hfl

\nhp It is easy to check that \se{\tet}{F}\ is a equivalence
relation in \q M. We will write

\nhp \hfl \q M/\cF\ \, = \, \{\q x/\cF\ : \q x \en\ \q M\}
 \hfl

\nhp the set of all equivalence classes of \se{\tet}{F} (\q M/\cF\
$\neq \empty$). If \ol x \en\ \po n M, we define

\nhp \hfl \ol x/\cF\ \, = \, \fm{x_1/{\cF}, \ldots, x_n/{\cF}}\ \,
\en\ \, \po{n}{(M/{\cF})}. \hfl

\nhp For each \ol x \en\ \po n M and \q i \en\ \q I, we take

\nhp \hfl \ol x(\q i) \, = \, \fm{x_1(i), \ldots, x_n(i)}\ \, \en\
\, \po{n}{M_i}. \hfl

\nhp With the notation in \ref{Lprod-ct}, if \q R \en\ $rel(n,
L)$, \op\ \en\ $op(n, L)$ and \ol x, \ol y \en\ \po{n}{M} are such
that
 \ol x/\cF\ = \ol y/\cF, then:

\nhp (A) \, \po{M}{\op}(\ol x)/\cF\ = \po{M}{\op}(\ol y)/\cF;

\nhp (B) \, \{\q i \en\ \q I : \se M i \mo\ $R[\ol x(i)]$\} \en\
\cF\ \hem \Lra\ \hem \{\q i \en\ \q I : \se M i \mo\ $R[\ol
y(i)]$\} \en\ \cF.

\nhp With the aid of (A) e (B) we can make \q M/\cF\ a \q
L-structure through the followings conditions:

\nhp $\ast$ If \q c \en\ $ct(L)$, \, \po{M/F}{c} \, = \,
\fm{\po{M_i}{c}}/\cF, i. e., the interpretation of the constant
symbol \q c in \q M/\cF\ is the equivalence class of the \q
I-sequence whose coordinates are the interpretations of \q c in
each component \se M i;

\nhp $\ast$ If \op\ \en\ $op(n, L)$ and \ol x \en\ \po n M, then
\, \po{M/F}{\op}(\ol x/\cF) \, = \, \po{M}{\op}(\ol x)/\cF;

\nhp $\ast$ If \q R \en\ $rel(n, L)$ and \ol x \en\ \po n M,

\nhp \hfl \q M/\cF\ \mo\ \q R[\ol x/\cF] \hem \Lra\ \hem \{\q i
\en\ \q I : \se M i \mo\ $R[\ol x(i)]$\} \, \en\ \, \cF. \hfl

\nhp Induction on complexity gives

\nhp (C) If \ta(\npl v) is a term in \q L and \ol x \en\ \po n M,
\, \po{M/F}{\ta}(\ol x/\cF) \, = \, \po{M}{\ta}(\ol x)/\cF.

\nhp (D) If \fh(\npl v) is an atomic formula in \q L and \ol x
\en\ \po n M,

\nhp \hfl \q M/\cF\ \mo\ \fh[\ol x/\cF] \hem \Lra\ \hem \{\q i
\en\ \q I : \se M i \mo\ $\fh[\ol x(i)]$\} \, \en\ \, \cF. \hfl

\nhp (E) The natural map  \q x \en\ \q M $\longmapsto$ \q x/\cF\
\en\ \q M/\cF\ is a surjective \q L-homomorphism.

\medskip

The \q L-structure \q M/\cF\ is named {\em the reduced product} of
the family \{\se M i : \q i \en\ \q I\} by the filter \cF. If \cF\
is an \und{ultrafilter} in \q I, \q M/\cF\ is called {\em the
ultraproduct} of \{\se M i : \q i \en\ \q I\} by  the ultrafilter
\cF. When all the \q L-structures are the same, \se M i = \q N, \q
i \en\ \q I, the correspondent construction is called {\em reduced
power} and {\em ultrapower}, when \cF\ is an ultrafilter, it is
indicated \po I N/\cF.

\ect

The fundamental result concerning ultraproducts is the:

\bte\ \label{Los-te} {\rm (\L\'os's Theorem)} Let $I$ a non empty
set, \{\se M i : \q i \en\ \q I\} is a family of non empty \q
L-structures, \q M = $\prod_{i\in I}$ \se M i  and \cF\  an {\em
ultrafilter} in \q I . Then for all formula $\fh(\npl v)$ in \q L
and all \ol x \en\ \po n M

\nhp $($\L$)$ \hfl \q M/\cF\ \mo\ $\fh[\ol x/F]$ \hem \Lra\ \hem
\{\q i \en\ \q I : \se M i \mo\ $\fh[\ol x(i)]$\} \, \en\ \, \cF.
\hfl \label{Los-equiv}

\ete

\bdm See Theorem 4.1.9, page 217, in \cite{CK}\ or Theorem 5.2.1,
page 90, in \cite{BS}. \qdr

\bre \label{Los-re} We add that the equivalence in the \L\'os's
Theorem remains true for reduced products in general (\cF\ {\em is
a filter}) if we restrict ourselves to formulas $\fh(\npl v)$ that
are in $p.p.(L)$ or, most generically, to the formulas that are
generated from the atomic formulas by the usage of the conjunction
and both quantifiers.

If \q M is a \q L-structure, \q I is a set e \cF\ $\sub P(I)$ is a
a filter in \q I, then there is a canonical \q L-homomorphism,
{\em the
diagonal from \q M to \po I M/\cF} \\
\hfl \de\ : \q M \lra\ \po I M/\cF, where \de(\q a) \, =
\, \fm{\q a}/\cF, \hfl\\
for each \q a \en\ \q M  in the equivalence class of the constant
\q I-sequence of value \q a.

It follows from \L\'os's Theorem  that when \cF\ is a ultrafilter
in \q I then the diagonal morphism, \de\ : \q M \lra\ \po I M/\cF,
is a elementary embedding. Similarly, if \cF\ is just a (proper)
filter in \q I then the diagonal morphism, \de\ : \q M \lra\ \po I
M/\cF, is a just a pure embedding (item \ref{Lpur-ssec}).

Another important consequence of this Theorem is that any
elementary class of structures is closed under the ultraproduct
construction.

\ere

\subsection{Colimits in L-mod} \label{Lcolim-ssec}

\bct{{\bf Filtered Colimits in  L-mod:}} \label{L-colim-ct} Let
\fm{I, \m}\ be a directed poset and {\cM} an \q I-diagram
 \footnote{As usual, we consider here \q I as a category
whose objects are the members of the set \q I and whose arrows are
the elements of the binary relation \m.}
in {\bf L-mod}.\\
\hfl $\cM : I \ \lra \ {\bf L-mod}$ \ : \ $(i \m j) \ \mapsto \
(M_{i} \overset{f_{ij}}\to\lra M_{j})$ \hfl

 Let \q W = $\cpd_{i \in I}$ \se{M}{i} = $\bcup_{i \in I}$
\xy{\se{M}{i}}{\{i\}}\ , be the disjunct reunion of the sets
\se{M}{i}. We have the canonical functions
 \se{w}{i}\ : \se{M}{i} \lra\ \q W, \q x $\mapsto$ \fm{x, i}. As
\q I is a directed poset the prescription

\nhp \hfl \fm{x, i}\ \equ\ \fm{y, j}\ \, \Lra\ \, \ex\ \q k \G\
$i, j$ \ such that \se{f}{ik}(\q x) = \se{f}{jk}(\q y), \hfl

\nhp defines an equivalence relation \equ\ in \q W. Let

\nhp \hfl \q M = \{\fm{x, i}/\equ\ \ : \ \fm{x, i}\ \en\ \q W\}
\hfl

\nhp be the set of all equivalence classes of \equ. Notice that
for each constant symbol \q c in \q L we have \fm{\po{{M}_{i}}{c},
i}\ \equ\ \fm{\po{{M}_{j}}{c}, j}. We interpret \q L in \q M as
follows: for each \q n \G\ 1 \ and \ \ol x \en\ \po{n}{M}, \ \ol x
= \fm{\fm{\se{x}{1}, \se{i}{1}}/\equ, \ldots, \fm{\se{x}{n},
\se{i}{n}}/\equ}, \ we define:

\nhp (A) If \q R \en\ $rel(n, L)$  then \q M \mo\ $R[\ol x]$ \,
iff

\nhp \hfl  \ex\ \q k \G\ \se{i}{1}, \ldots, \se{i}{n}, such that
\se{M}{k} \mo\ $R[\se{f}{\se{i}{1} k}(\se{x}{1}), \ldots,
\se{f}{\se{i}{n} k}(\se{x}{n})]$. \hfl

\nhp (B) If \op\ \en\ $op(n, L)$ we take \q k \G\ \se{i}{1},
\ldots, \se{i}{n}\ and define \po{M}{\op}(\ol x) as the
equivalence class of the pair

\nhp \hfl \fm{\po{{M}_{k}}{\op}(\se{f}{\se{i}{1} k}(\se{x}{1}),
\ldots, \se{f}{\se{i}{n} k}(\se{x}{n})), \, k}. \hfl

\nhp (C) If  \q c \en\ $ct(L)$ we take \, \po{M}{c}\ =
\fm{\po{{M}_{i}}{c}, i}/\equ.

\nhp Because \q I is directed, the constructions above are
independents of the particular chose of representations and also
of the index chose made above. Further, the compositions of the
quotient function, $q:  \q W \lra\ \q M$, with the functions
\se{w}{i}, defines {\em \q L-homomorphisms}  \se{\al}{i}\  :
\se{M}{i} \lra\ \q M \ that make (\q M, \{\se{\al}{i}\ : \q i \en\
\q I\}) a co-cone over the diagram {\cM} . This co-cone is the
colimit \ili\ \cM.

\ect

 \bco\ \label{Lcolim-co} Let {\cM}  =
$(M_i, \; \{\se{f}{ij}\ :\ i\ \m\ j\})$ be an \q I-diagram in {\bf
L-mod}, where \q I is a directed poset.
 A co-cone in {\bf L-mod} over
{\cM} , $(\q N,  \ \{\beta_{i} : i \in I\} )$, is
$(\text{isomorphic to})$ \ili\ {\cM} iff it verifies the following
conditions:

\nhp $[colim\ 1]$ : \q N = \bcup\ \{$\se{\bt}{i}(M_i)$ : \q i \en\
\q I\}.

\nhp $[colim\ 2]$ : If $\fh(\npl{v})$ is an atomic formula in \q L
and \ol s \en\ \po{n}{N},

\nhp \hfl  \q N \mo\ $\fh[\ol s]$ \hmm \Lra\ \hmm $\left\{
\begin{array}{c} \mbox{\ex\ \q k \en\ \q I \ and \ \ol x \en\
\po{n}{M_k}\
such that}\vsm\\
\mbox{\se{s}{p}\ = $\se{\bt}{k}(\se{x}{p})$,  \q 1 \m\ \q p \m\ \q n,} \vsm\\
 \mbox{and \,  $M_k$ \mo\ $\fh[ \ol x ]$.} \end{array} \right . $ \hfl

\qdr \eco

\subsection{Reduced products and filtered colimits of products}
\label{colu-ssec}

\hp There is a connection between reduced products and certain
filtered colimits\footnote{The geometrical girth of this result
appears in \cite{Ell} and \cite{Mir}.} that will be very useful in
the proof of our main result, namely Theorem \ref{repro-te}.
Before the precise statement and its proof we need to establish
some notation.

\bct\ \label{colu-ct} Let \q L be a first-order language with
equality, $I$ a non-empty set, \{\se M i : \q i \en\ \q I\} a
family of \q L-structures all non-empty and \q M = $\prd_{i \in
I}$ \se M i their product (item \ref{Lprod-ct}).

\nhp (A) For each \q J \sub\ \q I \, let \rs{M}{J}\ \, = \,
$\prd_{j \in J}$ \se M j;

\nhp (B) If \q J \sub\ \q K \sub\ \q I \, then there is a
canonical \q L-morphism, \se{\pi}{KJ}\ : \rs{M}{K}\ \lra\
\rs{M}{J}, that forgets the coordinates out of \q J, that is, for
\q x \en\ \rs{M}{K}, \se{\pi}{KJ}(\q x) = \rs{x}{J}\ (we recall
that \q x is a function from \q K to $\bcup_{k \in K}$ \se M k).
Regard that

\nhp (*) \hfl \se{\pi}{JJ}\ \, = \, $Id_{M_{|J}}$ \hmm and \hmm \q
J \sub\ \q K \sub\ \q W \, \Ra\ \, \se{\pi}{WJ}\ \, = \,
\se{\pi}{KJ}\ \com\ \se{\pi}{WK}. \hfl

\nhp The canonical projections, \se{\pi}{i}\ : \q M \lra\ \se M i,
correspond to \se{\pi}{I\{i\}}.

\nhp (C) For each \q J \sub\ \q I \, we define $\ast$ :
\xy{\rs{M}{J}}{M}\  \lra\ \q M, \fm{s, x}\ $\longmapsto$ \q s
$\ast$ \q x, where

\nhp \hfl \q s $\ast$ \q x(\q i) \, = \, $\left \{
\begin{array}{ll}
s(i) & \mbox{if \q i \en\ \q J}\vsm\\
x(i) & \mbox{if \q i \nen\ \q J.} \end{array} \right . $ \hfl

\nhp Note that when \q J = \q I then the operation $\ast$ is the
projection in the first coordinate. Equivalently, for each \q x
\en\ \q M, the function \
 (\dt) $\ast$ \q x : \q M  \lra\ \q M \ is the identity function.

Let \cF\ be a filter in \q I. Then \fm{\cF, \sub}\ is a downward
directed poset (item \ref{conjdir-ssec}) because for each \q J, \q
K \en\ \cF\ then  \q J \Int\ \q K \en\ \cF. Consequently,
\po{op}{\cF}, the opposite poset de \fm{\cF, \sub}\ , is (upward)
directed. Consider

\nhp \hfl \cM\ \, = \, (\rs{M}{J} ,  \{\se{\pi}{KJ}\ : \q J \sub\
\q K \, and \, \q J \en\ \cF\}). \hfl

\nhp By (*) in (B), \cM\ is a \po{op}{\cF}-diagram in {\bf L-mod},
the {\em directed diagram associated} to the family \{$M_i$: \q i
\en\ \q I\} and to the filter \cF\ in \q I.

\ect

\bpr \label{colu-pr} Let \q L be a first-order language with
equality, $I$ a non-empty set, \{\se M i : \q i \en\ \q I\} a
family of \q L-structures all non-empty and \cF\ a filter in \q I.
Consider

\nhp \hfl \cM\ \, = \, $(\rs{M}{J}$ ,  $\{\se{\pi}{KJ}\ : \q J
\sub\ \q K$ , $\q J \en\ \cF\})$ \hfl

\nhp the  directed diagram associated to the family \{\se M i : \q
i \en\ \q I\} and to the filter \cF\ in \q I as in
$\ref{colu-ct}$. Then \ili\ \cM\ is naturally \q L-isomorphic to
the reduced product $\prd_{i \in I}$ \se M i/\cF.

\epr

\bdm (Sketch) We make use of the notation in \ref{colu-ct}. In
particular, \q M \rs{M}{I}\ =  $\prd_{i \in I}$ \se M i \ is the
\q L-structure product.

\nhp (A) We fix a \q t \en\ \q M. For each \q J \en \cF, we
consider the mapping  \se{\nu}{J}\ : \rs M J \lra\ \q M/\cF, given
by

\nhp \hfl \se{\nu}{J}(\q s) \, = \, (\q s $\ast$ \q t)/\cF. \hfl

\nhp  As \und{$M \neq \empty$} this definition make sense. Follows
directly from the definition of reduced product that the function
\se{\nu}{J} \ {\em does not depends} of the particular element $t
\in M$ chosen.

\nhp (B) It may be verified, from the constructions of the product
structure, reduced product and filtered colimit that, for each \q
J \en\ \cF\ the function \se{\nu}{J}\ : \rs M J \lra\ \q M/\cF\ is
a \q L-homomorphism. Further, for each \q J, \q K \en\ \cF\ such
that \q J \sub\ \q K, the following diagram commutes:

\nhp \hfl \etri{\nf\nf \rs{M}{K}}{\rs{M}{J}}{\nf \q M/\cF}
{\put(-10,0){\se{\pi}{KJ}}}{\se{\nu}{J}}{\nf \se{\nu}{K}}\ \hfl

\nhp (C) It may be checked that   (\q M/\cF, \{\se{\nu}{J}\ : \q J
\en\ \cF\}) is a co-cone over the diagram \cM\ that
  satisfies a universal property and then it must be (isomorphic to) the
 co-cone \ili\ \cM, the colimit of the diagram \cM.
\qdr

\bre\ \label{colu-re} We note that if \ \cM\ \, = \, (\rs{M}{J},
\; \{\se{\pi}{KJ}\ : \q J \sub\ \q K , \q J \en\ \cF\}) \ is the
directed diagram of the proposition below then the $L$-structure \
$\li {\cM}$ \ seems to be the ``fundamental'' notion of reduced
product (or ultraproduct, when \cF\ is a ultrafilter) because this
is the structure that {\em always} is defined and that {\em
always} satisfies \L\'os's equivalence (\L), page
\pageref{Los-equiv}, and its version for reduced products (see
Theorem \ref{Los-te} and Remark \ref{Los-re}) . However, if we
admit just an empty structure $M_{j}$ in the original definition
of reduced product, take some filter {\cF} such that $\{ j \}
\notin {\cF}$ and regard the p.p.-sentence ``I am not the empty
structure'' : $\exists v_{0} (v_{0} = v_{0})$ , then we have that
\ $\prd_{i \in I}$ \se M i/\cF\ \, is empty \ but \ $\{ i \in I :
M_{i}$ is empty $\} \notin {\cF}$.

\ere

\section{Profinite Structures and  Ultraproducts}

\hp We present now ours results.

\bdf\ \label{sprofindef} A \q L-structure is {\bf profinite} when
it is \q L-isomorphic to the limit of a diagram of {\em finite} \q
L-structures over a {\em downward directed poset}.

\edf

\bre \label{reprofinre} If \q P is a profinite \q L-structure then
there is an {\em upward} directed poset, \fm{I, \m}, and a
cofiltered diagram of finite \q L-structures over \q I,

\nhp \hfl \cM\ \, = \, $(\se M i ,  \text{\{\se{f}{ji}\ : \q i \m\
\q j\}})$ \hfl

\nhp such that $(P , \{\lambda_{i} : \q i \en\ \q I\} )$ \, = \,
\dli\ \cM. By Proposition \ref{Llim-co} we can consider \q P as a
substructure of the product \q M = $\prd_{i \in I}$ \se M i, i.e.,
there is a natural \q L-imbedding, \ $\iota$ : \q P \lra\ \q M ,
such that for all \q i \en\ \q I, \,

\nhp $(\sharp)$ \hfl \se{\lam}{i}\ \, = \, \se{\pi}{i}\ \com\
$\iota$. \hfl \parbox{120pt}{\etri{\q P}{\q M}{\se M
i}{$\iota$}{\se{\pi}{i}} {\se{\lam}{i}}} \hfl

\nhp where \se{\pi}{i}\ : \q M \lra\ \se M i is the a canonical
projection. Further, it follows from $[lim\ 1]$ in \ref{Llim-co}
that

\nhp $(\flat)$ \hfl  \fa\ \ol x \en\ \q P \ \fa\ \q j,\q k \en\ \q
I \ ( \q j \en\ \up k \hem \Ra\ \hem \se{f}{jk}(\se x j) \, = \,
\se x k ). \hfl \label{flat}

We saw in \ref{colu-ct} that, if \cF\ is a filter in \q I then for
each \q J \en\ \cF\ we have a natural \q L-morphism

\nhp \hfl \se{\nu}{J}\ : \rs{M}{J}\ \lra\ \q M/\cF, given by \q x
$\longmapsto$ \q x/\cF, \hfl

\nhp where \q M/\cF\ indicates the reduced product $\prd_{i \in
I}$ \se M i/\cF.

\ere

With these preliminary we enunciate the

\bte \label{repro-te} Profinite \q L-structures are retracts of
ultraproducts of finite \q L-structures. More precisely, and with
the notation in $\ref{reprofinre}$, let \fm{I, \m}\ be a directed
poset and

\nhp \hfl \cM\ \, = \, $(\se M i , \{{f}_{ji}\ : \q i \m\ \q j\})$
\hfl

\nhp a cofiltered diagram of finite \q L-structures over \q I. If
\dli\ \cM\ \, = \, $(P , \{{\lambda}_{i} : \q i \en\ \q I\})$ then
the \q L-morphism that is the composition

\nhp \hfl \q P \, $\srel{\iota}{\lra}$ \, $\prd_{i \in I}$ \se M i
\, $\srel{\nu_I}{\lra}$ \, $\prd_{i \in I}$\se M i/\cU, \hfl

\nhp is an \q L-section (item \ref{retr-ssec}), where \cU\ is a
directed ultrafilter in \q I (item \ref{conjdir-ssec}) .

\begin{picture}(100,120)
\setlength{\unitlength}{.6\unitlength} \thicklines
\put(30,150){$P$} 
\put(120,150){$\prd_{i \in I} M_{i}$}
\put(250,150){$\prd_{i \in I} M_i/\cU$}
\put(260,15){$P$} 
\put(55,155){\vector(1,0){50}} \put(195,155){\vector(1,0){45}}
\put(275,130){\vector(0,-1){90}} \put(40,135){\vector(2,-1){210}}
\put(80,165){$\iota$} \put(210,165){$\nu_{I}$}
\put(285,70){$\gamma^{U}$}
\put(100,70){$Id_P$} %
\put(200,95){$\circlearrowright$} \put(400,75){Profinite and
ultraproduct of finite}
\end{picture}

\ete

\bdm  By Lemma \ref{ufildir-le} there is a directed ultrafilter in
\fm{I, \m}; the proof will be carried on fixing a such ultrafilter
\cU.

\nhp Let \q M = $\prd_{i \in I}$ \se M i be the product \q
L-structure of the family $\{\se M i : i \in I\}$. By the
Proposition \ref{colu-pr} (and with the same notation), we know
that

\nhp \hfl \q M/\cU \, = \, $\prd_{i \in I}$ \se M i/\cU\ \, is \q
L-isomorphic to \, \ili\ (\rs{M}{J}; \{\se{\pi}{KJ}\ : \q J \sub\
\q K , \q J \en\ \cU\}). \hfl

\nhp We shall use  this fact to build a \q L-morphism \po{U}{\ga}\
such that

\nhp \hfl \po{U}{\ga}\ \com\ (\se{\nu}{I} \com $\iota$) \,= \,
$Id_P$, \hfl

\nhp then the demonstration will be finished. As \cU\ will remain
fixed through the proof, we will indicate \po{U}{\gamma}\ just by
\ga. As the proof is a little bit long and technical it will be
carry through with the aid of several Facts. We will make free
usage of the notational conventions in \ref{Lprod-ct} and
\ref{reprofinre}.

\medskip

For each \q J \en\ \cU, \q i \en\ \q I, \ol x \en\ \rs{M}{J}\ =
$\prd_{j \in J}$ \se M j \ and \q y \en\ \se M i \, we define

\nhp \hfl \se{V}{J,i}(\ol x, \q y) \, = \, \{\q j \en\ \q J \Int\
\up i : \se{f}{ji}(\se x j) = \q y\}. \hfl

\bfa\ \label{reprofa1} For each \q J \en\ \cU, \q i \en\ \q I, \ol
x \en\ \rs{M}{J}\ \ and \q y, \q z \en\ \se M i,

\nhp a)  \, \q z \niq\ \q y \, \Ra\ \, $\se{V}{J,i}(\ol x, y)\
\Int\ \se{V}{J,i}(\ol x, z)$ \, = \, \0.

\nhp b) \, \q J \Int\ \up i \, = \, $\cpd_{y \in M_i}$
$\se{V}{J,i}(\ol x, y)$. \footnote{\cpd\ indicates that this union
is {\em disjunctive}.}

\efa

\nhp {\em Proof.}  Item (a) follows immediately from the fact that
\se{f}{ji}\ is a function. For (b), by the definition of
\se{V}{J,i}(\ol x, \q y) it is clearly enough to show that the
left side of the equality is contained in the right side, but note
that if \q j \en\ \q J \Int\ \up i \, then \se{f}{ji}(\se x j)
\en\ \se M i, as required. \qdr

\bfa\ \label{reprofa2} For each \q J \en\ \cU\ and \q i \en\ \q I
\ there is a \q L-morphism

\nhp \hfl \se{\ga}{J,i}\ : \rs{M}{J}\ = $\prd_{k \in J}$ \se M k
\lra\ \se M i \hfl

\nhp such that

\nhp a) If \ol x \en\ \rs{M}{J}\ and \q y \en\ \se M i \, then \,
$\se{\ga}{J,i}(\ol x)\ =\ y$ \hem iff \hem $\se{V}{J,i}(\ol x, y)$
\en\ \cU.

\nhp b) If \q J \sub\ \q K are members of \cU\ \ and \q i \en\ \q
I \ then the left diagram below commutes:

\nhp \hfl \parbox{120pt}{\etri{\nf \nf \rs{M}{K}}{\rs{M}{J}}{\se M
i} {\nf \se{\pi}{KJ}}{\se{\ga}{K,i}}{\nf\nf \se{\ga}{J,i}}} \hfl
\parbox{120pt}{\etri{\nf \rs{M}{J}}{\se M k}{\se M i}{\nf \se{\ga}{J,k}}
{\se{f}{ki}}{\nf \se{\ga}{J,i}}} \hfl

\nhp c) For each \q J \en\ \cU\ and \q i \m\ \q k in \q I \ the
right diagram below commutes.

\nhp d) For each \q k \en\ \q I, \, \se{\ga}{I,k}\ \com\ $\iota$
\, = \, \se{\pi}{k}\ \com\ $\iota$, where \se{\pi}{k}\ : \q M
\lra\ \se M k \ is the canonical projection.

\efa

\nhp {\em Proof.} Because \cU\ is a directed filter in \fm{I, \m}
(item \ref{conjdir-ssec}) for each  $J \in {\cU}$ and $ i \in I$
we have \q J \Int\ \up i \en\ \cU\ ; because \cU\ is an
ultrafilter and \se M i is finite, the Fact \ref{reprofa1}.(b)
implies that {\em there is a unique} \q y \en\ \se M i such that
\se{V}{J,i}(\ol x, \q y) \en\ \cU. We define

\nhp \hfl \se{\ga}{J,i}(\ol x) \, = \, the unique \q y \en\ \se M
i such that \se{V}{J,i}(\ol x, \q y) \en\ \cU. \hfl

\nhp It is clear that the item (a) is verified. Now, we must show
that \se{\ga}{J,i}\ is a \q L-morphism. To make easier the
reading, if \q J \en\ \cU, we will indicate the symbols
interpretations of \q L in \rs{M}{J}\ by an exponent \q J. Then,
if \q c is a constant symbol in \q L, we will use \po{J}{c}\
instead  \po{M_{|J}}{c}; analogously for the functional and
relational symbols.

\nhp $\ast$ Let \q c \en\ $ct(L)$. We saw in \ref{Lprod-ct} (item
\ref{Llim-ssec}) that \po J c is a sequence \fm{\po{M_j}{c}}\ \en\
\rs{M}{J}. So, as the \se{f}{ji}\ are \q L-morphisms, we get

\nhp \hfl \se{V}{J,i}(\po J c, \po{M_i}{c}) \, = \, \{\q j \en\ \q
J \Int\ \up i : \se{f}{ji}(\se{c^J}{j}) = \po{M_i}{c}\} \, = \,
\{\q j \en\ \q J \Int\ \up i : \se{f}{ji}(\po{M_j}{c}) =
\po{M_i}{c}\} \, = \, \q J \Int\ \up i \hfl

\nhp that belongs to \cU. By item (a), \se{\ga}{J,i}(\po J c) =
\po{M_i}{c}, as we wish.

\nhp $\ast$ Let \op\ be a \q n-ary functional symbol in \q L. If
\npl{\ol x}\ \en\ \po{n}{(\rs{M}{J})}\ and \q j \en\ \q J then, by
\ref{Lprod-ct}, we have

\nhp (A) \hfl \po{J}{\op}(\npl{\ol x})(\q j) \, = \,
\po{M_j}{\op}(\se{x}{1j}, \ldots, \se{x}{nj}). \hfl

\nhp Consider

\nhp \hfl $\left \{ \begin{array}{lcll}
\po p y & = & \se{\ga}{J,i}(\se{\ol x}{p}), & \mbox{1 \m\ \q p \m\ \q n;}\vsm\\
z & = & \po{M_i}{\op}(y^1, \ldots, y^n); & \vsm\\
h & = & \po{J}{\op}(\npl{\ol x}) & \mbox{(\en\ \rs M J).}
\end{array} \right . $ \hfl

\nhp We will show that

\nhp (B) \hfl $\bcap_{p = 1}^{n}$ \se{V}{J,i}(\se{\ol x}{p}, \po p
y) \, \sub\ \, \se{V}{J,i}(\q h, \q z). \hfl

\nhp If \q j \en\ $\bcap_{p = 1}^{n}$ \se{V}{J,i}(\se{\ol x}{p},
\po p y) \, then the definition of \se{V}{J,i}\ implies

\nhp (C) \hfl \fa\ 1 \m\ \q p \m\ \q n, \, \se{f}{ji}(\se{x}{pj})
\, = \, \po p y. \hfl

\nhp Because the \se{f}{ji}\ are \q L-morphisms, (A) and (C) give

\nhp $\se{f}{ji}(h_j)$ \, = \,
\se{f}{ji}(\po{M_j}{\op}(\se{x}{1j}, \ldots, \se{x}{nj})) \, = \,
\po{M_i}{\op}(\se{f}{ji}(\se{x}{1j}), \ldots,
\se{f}{ji}(\se{x}{nj})) \, = \, \po{M_i}{\op}(\po 1 y, \ldots, \po
n y) \, = \, \q z,

\nhp and this proves (B). As the intersection of the left side in
(B) belongs to \cU \, we have \se{V}{J,i}(\q h, \q z) \en\ \cU. By
the item (a) of this Fact, this means that

\nhp \hfl \se{\ga}{J,i}(\po{J}{\op}(\npl{\ol x})) \, = \,
\po{M_i}{\op}(\se{\ga}{J,i}(\se{\ol x}{1}, \ldots, \se{\ol x}{n}))
\hfl

\nhp showing that \se{\ga}{J,i}\ preserves the operation \op;

\nhp $\ast$ Let \q R  be a \q n-ary relational symbol in \q L.
Consider \npl{\ol x}\ \en\ \po{n}{(\rs{M}{J})}. By \ref{Lprod-ct}

\nhp (D) \hfl \rs{M}{J}\ \mo\ \q R[\npl{\ol x}] \, iff \, \fa\ \q
j \en\ \q J, \, \se M j \mo\ \q R[\se{x}{1j}, \ldots, \se{x}{nj}].
\hfl

\nhp As above, let \po p y = \se{\ga}{J,i}(\se{\ol x}{p}), 1 \m\
\q p \m\ \q n. We must show that

\nhp (E) \hfl \rs{M}{J}\ \mo\ \q R[\npl{\ol x}] \hem \Ra\ \hem \se
M i \mo\ \q R[\po 1 y, \ldots, \po n y]. \hfl

\nhp Because $\bcap_{p = 1}^{n}$ \se{V}{J,i}(\se{\ol x}{p}, \po p
y) \, \en\ \, \cU, this intersection is non-empty; if \q j is a
member of this intersection, the topic (C) above is checked. Then,
it follows from (D) and the fact that \se{f}{ji}\ is a \q
L-morphism that

\nhp \hfl \rs{M}{J}\ \mo\ \q R[\npl{\ol x}] \hem \Ra\ \hem \se M j
\mo\ \q R[\se{x}{1j}, \ldots, \se{x}{nj}] \hem \Ra\ \hem \se M i
\mo\ \q R[\se{f}{ji}(\se{x}{1j}), \ldots, \se{f}{ji}(\se{x}{nj})],
\hfl

\nhp with this and (C) we obtain (E), completing the proof that
\se{\ga}{J, i}\ is a \q L-morphism.

\nhp b) Let \ol t \en\ \rs{M}{K}\ and \ol x = \se{\pi}{KJ}(\ol t)
\footnote{We recall that \se{\pi}{KJ}\ is the projection that
forgets the coordinates out of \q K.}. If \q y = \se{\ga}{J,i}(\ol
x) \, we will see that

\nhp \hfl \se{V}{J,i}(\ol x, \q y) \, \sub\ \, \se{V}{K,i}(\ol t,
\q y). \hfl

\nhp In fact, if \q j \en\ \se{V}{J,i}(\ol x, \q y) (obviously
contained \q K \Int\ \up i) then

\nhp \hfl \se{f}{ji}(\se t j) \, = \, \se{f}{ji}(\se x j) \, = \,
\q y, \hfl

\nhp as required. As \se{V}{J,i}(\ol x, \q y) \en\ \cU \, we have
\se{V}{K,i}(\ol t, \q y) \en\ \cU\ and the item (a) ensures that
\se{\ga}{K,i}(\ol t) = \q y = \se{\ga}{J,i}(\se{\pi}{KJ}(\ol t)),
as we need.

\nhp c) Let \ol x \en\ \rs{M}{J}\ and \q z = \se{\ga}{J,k}(\ol x).
Then

\nhp (F) \hfl \se{V}{J,k}(\ol x, \q z) \, \sub\ \, \se{V}{J,i}(\ol
x, \se{f}{ki}(\q z)). \hfl

\nhp In fact, if \q j \en\ \se{V}{J,k}(\ol x, \q z) (contained in
\q J \Int\ \up i because \q i \m\ \q k) then \se{f}{jk}(\se x j) =
\q z. As \cM\ is a cofiltered diagram, we have

\nhp \hfl \se{f}{ji}(\se x j) \, = \, \se{f}{ki}(\se{f}{jk}(\se x
j)) \, = \, \se{f}{ki}(\q z) \hfl

\nhp showing that \q j \en\ \se{V}{J,i}(\ol x, \se{f}{ki}(\q z));
as the topic (F) above ensures that this set belongs to \q U, the
item (a) implies  \se{\ga}{J,i}\ \, = \, \se{f}{ki}\ \com\
\se{\ga}{J,k}, as needed.

\nhp d) For each \ol x \en\ \q P  and \q k \en\ \q I \, observe
that $\se{\pi}{k}(\iota(\ol x))$ = \se x k. It follows from the
relation $(\flat)$ in \ref{reprofinre} (page \pageref{flat}) that

\nhp \hfl \se{V}{I,k}($\iota(\ol x)$, \se x k) \, = \, \{\q j \en\
\up k : \se{f}{jk}(\se x j) = \se x k\} \, = \, \up k. \hfl

\nhp Because \cU\ is a directed ultrafilter, we have
\se{V}{I,k}($\iota(\ol x)$, \se x k) \en\ \cU\ and the item (a)
gives the needed conclusion, closing the proof of the Fact
\ref{reprofa2}. \qdr

\medskip

By Proposition \ref{colu-pr} we have

\nhp \hfl \q M/\cU\ \, = \, \ili\ $(\rs{M}{J} , \{ \se{\pi}{KJ} :
\q J \sub\ \q K , \q J \en\ \cU \})$ . \hfl

Fact \ref{reprofa2}.(b) and the universal property of the filtered
colimits ensures that, for each \q i \en\ \q I, {\em there is a
unique} \q L-morphism, \se{\ga}{i}\ : \q M/\cU \, \lra\ \, \se M
i, such that for all \q J \en\ \cU\ the left diagram below
commutes:

\nhp (*) \hfl \parbox{120pt}{\etri{\nf\nf\rs{M}{J}}{\q M/\cU}{\se
M i} {\se{\nu}{J}}{\se{\ga}{i}}{\nf \se{\ga}{J,i}}} \hfl
\parbox{120pt}{\etri{\nf \nf \q M/\cU}{\se M k}{\se M i}{\se{\ga}{k}}
{\se{f}{ki}}{\se{\ga}{i}}} \hfl \label{*}

\bfa\ \label{reprofa3} For each \q i \m\ \q k in \q I, the right
diagram above in $(*)$ is commutative.

\efa

\nhp {\em Proof.} For each \q i \m\ \q k in \q I and \q J \en\
\cU, the Fact \ref{reprofa2}.(c) gives \se{\ga}{J,i}\ \, = \,
\se{f}{ki}\ \com\ \se{\ga}{J,k}. Then, the commutativity of the
left diagram above in (*) $-$ for \q k and \q i $-$, implies that,
for all \q J \en\ \cU\, we have

\nhp \hfl \se{f}{ki}\ \com\ \se{\ga}{k}\ \com\ \se{\nu}{J}\ \, =
\, \se{f}{ki}\  \com\ \se{\ga}{J,k}\ \, = \, \se{\ga}{J,i}\ \, =
\, \se{\ga}{i}\ \com\ \se{\nu}{J}. \hfl

\nhp Now, the uniqueness of the \se{\ga}{i}\ that make the left
diagram commutative, for all \q J \en\ \cU\, ensures that
\se{f}{ki}\ \com\ \se{\ga}{k}\ = \se{\ga}{i}, as required. \qdr

\medskip

Fact \ref{reprofa3} shows that  $( M/\cU , \{\se{\ga}{i}\ : \q i
\en\ \q I \} )$ \ is a cone over the cofiltered diagram \cM. Then
the universal property of the cofiltered limits ensures that {\em
there is a unique} \q L-morphism

\nhp (**) \hfl \ga\ : \q M/\cU\ \, \lra\ \, \q P = \dli\ \cM \hfl
\parbox{120pt}{\etri{\nf\nf \q M/\cU}{\q P}{\se M i}{\ga}{\se{\lam}{i}}
{\se{\ga}{i}}}\hfl

\nhp such that for all \q i \en\ \q I the diagram above comutes.

We will check now that

\nhp (G) \hfl \ga\ \com\ \se{\nu}{I}\ \com\ $\iota$ \, = \,
$Id_P$. \hfl

\nhp As $( P , \{\se{\lam}{i}\ : \q i \en\ \q I \})$  \, = \,
\dli\ \cM, the universal property of the limits ensures that to
prove (G) it is enough to show that for all \q k \en\ \q I

\nhp (H) \hfl \se{\lam}{k}\ \com\ (\ga\ \com\ \se{\nu}{I}\ \com\
$\iota$) \, = \, \se{\lam}{k}. \hfl

\nhp As \, \hfl $\left \{ \begin{array}{lcll} \se{\lam}{k}\ \com\
\ga\ & = & \se{\ga}{k}  & \mbox{by the diagram in (**);}
\vsm\\
\se{\ga}{k}\ \com\ \se{\nu}{I}\ & = & \se{\ga}{I,k}  &
\mbox{by the left diagram in (*), page \pageref{*};}\vsm\\
\se{\pi}{k}\ \com\ \iota\ & = & \se{\lam}{k}  &\mbox{by $(\sharp)$
in \ref{reprofinre}, page \pageref{flat},}
\end{array} \right . $ \hfl

\nhp  (H) is equivalent to \se{\ga}{I,k}\ \com\ $\iota$ =
\se{\pi}{k}\ \com\ $\iota$, but that is precisely the content of
the Fact \ref{reprofa2}.(d), so the proof of the Theorem is
complete. \qdr

\nhp \hfl \begin{picture}(200,200)
\setlength{\unitlength}{.6\unitlength} \thicklines
\put(25,20){$M_{k}$}
\put(25,280){\q P} 
\put(150,280){\q M} 
\put(265,280){\q M/\cU} 
\put(265,20){\q P} 
\put(50,285){\vector(1,0){80}} \put(80,295){$\iota$}
\put(180,285){\vector(1,0){75}} \put(205,295){$\nu_{I}$}
\put(30,260){\vector(0,-1){200}} \put(00,150){\se{\lam}{k}}
\put(145,260){\vector(-1,-2){100}} \put(65,170){$\gamma_{I,k}$}
\put(255,260){\vector(-1,-1){200}} \put(185,150){$\gamma_{k}$}
\put(275,260){\vector(0,-1){200}} \put(290,150){\ga}
\put(250,25){\vector(-1,0){180}} \put(150,00){\se{\lam}{k}}
\put(75,200){$\circlearrowright$}
\put(160,200){$\circlearrowright$}
\put(180,70){$\circlearrowright$}
\end{picture} \hfl



\vst The Lemma \ref{Lpur-le} and the Theorem \ref{repro-te}
produce the following

\bco \label{repro-co} Let \cA\ = ${\rm Mod}(T)$ where $T$ is a
theory axiomatized by $L$-sentences of the form $\all \vec{x}
(\psi_{0}(\vec{x}) \ra \psi_{1}(\vec{x}))$ where
$\psi_{0}(\vec{x}),\psi_{1}(\vec{x})$ are formulas in
$\exists^{+}(L)$. Then the full subcategory \cA\ $\sub$\ {\bf
L-mod} has profinite objects \footnote{We note that {\cA} has
already some finite object: \fo\ , the final object of {\bf
L-mod}, belongs to {\cA}, as we have noted in \ref{Lprod-ct} .},
that is, \cA\ is closed in {\bf L-mod} under the formation of such
limits. \qdr \eco

\section{Final Remarks} \label{fin-sec}

\hp As an application of the results above we mention the case of
the {\em Special Groups}, a first-order axiomatization of the
algebraic theory of quadratic forms (see \cite{DM2}). The suitable
first-order language, $L_{SG}$, contains two symbols for constants
(1 and -1), one symbol for binary operation (multiplication) and
one symbol for quaternary relation ($\equiv$, the isometry between
quadratic forms with dimension 2). The special groups axioms
(Definition 1.2 in \cite{DM2}) have the form $\all \vec{x}
(\psi_{0}(\vec{x}) \ra \psi_{1}(\vec{x}))$, where
$\psi_{0}(\vec{x}),\psi_{1}(\vec{x})$ are existencial-positives
$L_{SG}$-formulas, from the results \ref{repro-te} and
\ref{repro-co} above we can conclude that there are profinite
special groups and that they are retracts of ultraproducts of
finite special groups, a result contained in Theorem 5.8 in
\cite{Mrn1} and that has further consequences (e.g., in the
forthcoming \cite{MDM}).

As we mentioned before, our main motivation in \cite{Mrn1} was to
study the class of profinite special groups and particularly the
construction of the Profinite Hull of Special Groups functor (an
English language version of those main results will appear in
\cite{MM}; see also the forthcoming \cite{Mrn2} for further
application of such construction). The perception that some of
those constructions and results can be transported toward the
general context of $L$-structures has appeared with the results of
the present work. Further material about the ''Model Theory of
Profinite Structures'' are being elaborated in \cite{Mrn3}.

\section*{Acknowledgements}

\hp I am in great debt with Prof. Francisco Miraglia of IME-USP
(Instituto de Matem\'atica e Estat\y stica, Universidade de S\~ao
Paulo, Brazil) who gave me many helpful suggestions to improve the
exposition of this work. I also would like to thank Prof. Walter
A. Carnielli of CLE-Unicamp (Centro de L\'ogica, Epistemologia e
Hist\'oria da Ci\^encia, Universidade de Campinas, Brazil) for his
generous help. I am also grateful to the Institutes IME-USP, were
I developed my PhD thesis, and CLE-Unicamp, were I now have the
opportunity of developing a post-doctoral fellowship with the
financial support of FAPESP (Funda\cao\ de Amparo \`a Pesquisa do
Estado de S\ao\ Paulo, Brazil).

\end{document}